\documentclass[12pt, reqno]{amsart}

\newtheorem{theorem}{Theorem}
\newtheorem{lemma}[theorem]{Lemma}

\usepackage{fullpage}

\newcommand{\sfrac}[2]{{\textstyle \frac {#1}{#2}}}

\title[A binary additive equation]{A binary additive equation involving fractional powers}
\author{Angel V. Kumchev}
\address{Department of Mathematics, Towson University, Towson, MD 21252-0001, U.S.A.}
\email{akumchev@towson.edu}

\date{May 5, 2006}

\begin{document}

\maketitle

\section{Introduction}

It is well-known that the number of integers $n \le x$ that can be expressed as sums of two squares is $O\big( x(\log x)^{-1/2} \big)$.
On the other hand, Deshouillers \cite{Desh73} showed that when $1 < c < \frac 43$, every sufficiently large integer $n$ can be represented
in the form
\begin{equation}\label{i.1}
  [m_1^c] + [m_2^c] = n,
\end{equation}
with integers $m_1, m_2$; henceforth, $[\theta]$ denotes the integral part of $\theta$. Subsequently, the range for $c$ in this result
was extended by Gritsenko \cite{Grit92} and Konyagin \cite{Kony03}. In particular, the latter author showed that \eqref{i.1} has solutions
in integers $m_1, m_2$ for $1 < c < \frac 32$ and $n$ sufficiently large.

The analogous problem with prime variables is considerably more difficult, possibly at least as difficult as the binary Goldbach problem.
The only progress in that direction is a result of Laporta \cite{Lapo99a}, which states that if $1 < c < \sfrac {17}{16}$, then almost all $n$
(in the sense usually used in analytic number theory) can be represented in the form \eqref{i.1} with primes $m_1, m_2$. Recently, Balanzario,
Garaev and Zuazua \cite{BaGaZuX} considered the equation
\begin{equation}\label{i.2}
  [m^c] + [p^c] = n,
\end{equation}
where $p$ is a prime number and $m$ is an integer. They showed that when $1 < c < \sfrac {17}{11}$, this hybrid problem can be solved for
almost all $n$. It should be noted that in regard to the range of $c$, this result goes even beyond Konyagin's. On the other hand, when $c$
is close to $1$, one may hope to solve \eqref{i.2} for all sufficiently large $n$, since the problem is trivial when $c = 1$. The main purpose of the present note is to address this issue. We establish the following theorem.

\begin{theorem}\label{th1}
  Suppose that $1 < c < \frac {16}{15}$. Then every sufficiently large integer $n$ can be represented in the form \eqref{i.2}.
\end{theorem}

The main new idea in the proof of this theorem is to translate the additive equation \eqref{i.2} into a problem about 
Diophantine approximation. The same idea enables us to give also a simple proof of a slightly weaker version of the result of Balanzario, Garaev and Zuazua. For $x \ge 2$, let $E_c(x)$ denote the number of integers $n \le x$ that cannot be represented in the form \eqref{i.2}. We prove the following theorem.

\begin{theorem}\label{th2}
  Suppose that $1 < c < \frac 32$ and $\varepsilon > 0$. Then
  \[
    E_c(x) \ll x^{3(1 - 1/c) + \varepsilon}.
  \]
\end{theorem}

We remark that Theorem \ref{th1} is hardly best possible. It is likely that more sophisticated exponential sum estimates and/or sieve techniques would have allowed us to extend the range of $c$. The resulting improvement, however, would have been minuscule; thus, we decided not to pursue such ideas.

\medskip

\noindent
\textbf{Notation.} Most of our notation is standard. We use Landau's $O$-notation, Vinogradov's $\ll$-symbol, and
occasionally, we write $A \asymp B$ instead of $A \ll B \ll A$. We also write $\{ \theta \}$ for the fractional part
of $\theta$ and $\| \theta \|$ for the distance from $\theta$ to the nearest integer. Finally, we define $e(\theta)
= \exp(2\pi i\theta)$.

\section{Proof of Theorem \ref{th1}: initial stage}
\label{s2}

In this section, we only assume that $1 < c < 2$. We write $\gamma = 1/c$ and set
\begin{equation}\label{2.1}
  X = \big( \sfrac 12 n \big)^{\gamma}, \quad X_1 = \sfrac 54X, \quad \delta = \gamma X^{1 - c}.
\end{equation}
If $n$ is sufficiently large, it has at most one representation of the form \eqref{i.2} with $X < p \le X_1$. Furthermore, such a
representation exists if and only if there is an integer $m$ satisfying the inequality
\begin{equation}\label{2.2}
  \big( n - [p^c] \big)^{\gamma} \le m < \big( n + 1 - [p^c] \big)^{\gamma}.
\end{equation}
We now proceed to show that such an integer exists, if $p$ satisfies the conditions
\begin{equation}\label{2.3}
  X < p \le X_1, \quad \{ p^c \} < \sfrac 12, \quad 1 - \sfrac 56\delta < \big\{ \big( n - p^c \big)^{\gamma} \big\} < 1 - \sfrac 23\delta.
\end{equation}
Under these assumptions, one has
\[
  X^{1 - c} = (n - X^c)^{\gamma - 1} < \big( n - p^c \big)^{\gamma - 1} \le (n - X_1^c)^{\gamma - 1} < 1.1X^{1 - c}.
\]
Hence,
\begin{align*}
  \big( n - [p^c] \big)^{\gamma}
  &= \big( n - p^c \big)^{\gamma} \left( 1 + \gamma\{ p^c \}\big( n - p^c \big)^{-1} + O \big( n^{-2} \big) \right) \\
  &< \big( n - p^c \big)^{\gamma} + \sfrac 12\gamma\big( n - p^c \big)^{\gamma - 1} + O \big( n^{\gamma - 2} \big) \\
  &< \big( n - p^c \big)^{\gamma} + 0.55\delta + O \big( \delta n^{-1} \big) \\
  &< \big[ \big( n - p^c \big)^{\gamma} \big] + 1 - 0.1\delta,
\end{align*}
and
\begin{align*}
  \big( n + 1 - [p^c] \big)^{\gamma}
  &= \big( n - p^c \big)^{\gamma} \left( 1 + \gamma(1 + \{ p^c \}) \big( n - p^c \big)^{-1} + O \big( n^{-2} \big) \right) \\
  &\ge \big( n - p^c \big)^{\gamma} + \gamma\big( n - p^c \big)^{\gamma - 1} + O \big( n^{\gamma - 2} \big) \\
  &> \big( n - p^c \big)^{\gamma} + \delta + O \big( \delta n^{-1} \big) \\
  &> \big[ \big( n - p^c \big)^{\gamma} \big] + 1 + 0.1\delta.
\end{align*}
Consequently, conditions \eqref{2.3} are indeed sufficient for the existence of an integer $m$ satisfying \eqref{2.2}. It remains to
show that there exist primes satisfying the inequalities in \eqref{2.3}. To this end, it suffices to show that
\begin{equation}\label{2.4}
  \sum_{X < p \le X_1} \Phi\big( p^c \big) \Psi \big( (n - p^c)^{\gamma} \big) > 0
\end{equation}
for some smooth, non-negative, $1$-periodic functions $\Phi$ and $\Psi$ such that $\Phi$ is supported in $(0, 1/2)$ and $\Psi$ is
supported in $(1 - \frac 56\delta, 1 - \frac 23\delta)$.

Let $\psi_0$ be a non-negative $C^{\infty}$-function that is supported in $[0, 1]$ and is normalized in $L^1$: $\| \psi_0 \|_1 = 1$.
We choose $\Phi$ and $\Psi$ to be the 1-periodic extensions of the functions
\[
  \Phi_0(t) = \psi_0(2t) \qquad \text{and} \qquad \Psi_0(t) = \psi_0\big( 6\delta^{-1}(t - 1) + 5),
\]
respectively. Writing $\hat\Phi(m)$ and $\hat\Psi(m)$ for the $m$th Fourier coefficients of $\Phi$ and $\Psi$, we can report that
\begin{equation}\label{2.6}
\begin{split}
  \hat\Phi(0) = \sfrac 12, \quad &|\hat\Phi(m)| \ll_r (1 + |m|)^{-r} \quad \text{for all } r \in \mathbb Z, \\
  \hat\Psi(0) = \sfrac 16\delta, \quad &|\hat\Psi(m)| \ll_r \delta(1 + \delta|m|)^{-r} \quad \text{for all } r \in \mathbb Z.
\end{split}
\end{equation}
Replacing $\Phi(p^c)$ and $\Psi ( (n - p^c)^{\gamma} )$ on the left side of \eqref{2.4} by their Fourier expansions, we obtain
\begin{equation}\label{2.7}
  \sum_{X < p \le X_1} \Phi\big( p^c \big) \Psi \big( (n - p^c)^{\gamma} \big) =
  \sum_{h \in \mathbb Z} \sum_{j \in \mathbb Z} \sum_{X < p \le X_1} \hat\Phi(h) \hat\Psi(j) \, e \big( hp^c + j(n - p^c)^{\gamma} \big).
\end{equation}
Set $H = X^{\varepsilon}$ and $J = X^{c - 1 + \varepsilon}$, where $\varepsilon > 0$ is fixed.
By \eqref{2.6} with $r = [\varepsilon^{-1}] + 2$, the contribution to the the right side of \eqref{2.7} from the terms with
$|h| > H$ or $|j| > J$ is bounded above by a constant depending on $\varepsilon$. Thus,
\[
  \sum_{X < p \le X_1} \Phi\big( p^c \big) \Psi \big( (n - p^c)^{\gamma} \big) =
  \sfrac 1{12}\delta \big( \pi(X_1) - \pi(X) \big) + O \big( \delta\mathcal R + 1 \big),
\]
where $\pi(X)$ is the number of primes $\le X$ and
\[
  \mathcal R = \mathop{\sum_{|h| \le H} \sum_{|j| \le J}}_{(h, j) \ne (0, 0)}
  \bigg| \sum_{X < p \le X_1} e \big( hp^c + j(n - p^c)^{\gamma} \big) \bigg|.
\]
Thus, it suffices to show that
\begin{equation}\label{2.5}
  \sum_{X < p \le X_1} e \big( hp^c + j(n - p^c)^{\gamma} \big) \ll X^{2 - c - 3\varepsilon}
\end{equation}
for all pairs of integers $(h, j)$ such that $|h| \le H$, $|j| \le J$, and $(h, j) \ne (0, 0)$.

\section{Bounds on exponential sums}

In this section, we establish estimates for bilinear exponential sums, which we shall need in the proof of \eqref{2.5}.
Our first lemma is a variant of van der Corput's third-derivative estimate (see \cite[Corollary 8.19]{IwKo04}).

\begin{lemma}\label{l1}
  Suppose that $2 \le F \le N^{3/2}$, $N < N_1 \le 2N$, and $0 < \delta < 1$. Let $f \in C^3[N, N_1]$ and suppose that we can
  partition $[N, N_1]$ into $O(1)$ subintervals so that on each subinterval one of the following sets of conditions holds:
  \begin{enumerate}
    \item  [i)] $\delta FN^{-2} \ll |f''(t)| \ll FN^{-2}$;
    \item [ii)] $\delta FN^{-3} \ll |f'''(t)| \ll FN^{-3}$, $|f''(t)| \ll \delta FN^{-2}$.
  \end{enumerate}
  Then
  \[
    \sum_{N < n \le N_1} e (f(n)) \ll \delta^{-1/2} \big( F^{1/6}N^{1/2} + F^{-1/3}N \big).
  \]
\end{lemma}

\begin{proof}
  Let $\eta$ be a parameter to be chosen later so that $0 < \eta \le \delta$ and let $\mathbf I$ be one of the subintervals
  of $[N, N_1]$ mentioned in the hypotheses. If i) holds in $\mathbf I$, then by \cite[Corollary 8.13]{IwKo04},
  \begin{equation}\label{3.3}
    \sum_{n \in \mathbf I} e(f(n)) \ll \delta^{-1/2}\big( F^{1/2} + NF^{-1/2} \big).
  \end{equation}

  Now suppose that ii) holds in $\mathbf I$. We subdivide $\mathbf I$ into two subsets:
  \[
    \mathbf I_1 = \big\{ t \in \mathbf I : \eta FN^{-2} \le |f''(t)| \ll \delta FN^{-2} \big\}, \quad
    \mathbf I_2 = \mathbf I \setminus \mathbf I_1.
  \]
  Since $f''$ is monotone on $\mathbf I$, the set $\mathbf I_1$ consists of at most two intervals and $\mathbf I_2$ is a
  (possibly empty) subinterval of $\mathbf I$. If $\mathbf I_2 = [a, b]$, then there is a $\xi \in (a, b)$ such that
  \[
    f''(b) - f''(a) = (b - a)f'''(\xi) \qquad \implies \qquad b - a \ll \eta \delta^{-1}N.
  \]
  Thus, by \cite[Corollary 8.13]{IwKo04} and \cite[Corollary 8.19]{IwKo04},
  \begin{gather}
    \sum_{n \in \mathbf I_1} e(f(n)) \ll \eta^{-1/2}\big( F^{1/2} + NF^{-1/2} \big),  \label{3.4a}\\
    \sum_{n \in \mathbf I_2} e(f(n)) \ll \eta\delta^{-4/3}F^{1/6}N^{1/2} + \eta^{1/2}\delta^{-2/3}F^{-1/6}N. \label{3.4b}
  \end{gather}
  Combining \eqref{3.3}--\eqref{3.4b}, we get
  \begin{equation}\label{3.5}
    \sum_{N < n \le N_1} e(f(n)) \ll \eta^{-1/2}\big( F^{1/2} + NF^{-1/2} \big) + \eta\delta^{-4/3}N^{1/2}F^{1/6}
    + \eta^{1/2}\delta^{-2/3}NF^{-1/6}.
  \end{equation}
  We now choose
  \[
    \eta = \delta \max\big( F^{-1/3}, F^{2/3}N^{-1} \big).
  \]
  With this choice, \eqref{3.5} yields
  \[
    \sum_{N < n \le N_1} e(f(n)) \ll \delta^{-1/2} \big( F^{1/6}N^{1/2} + F^{-1/3}N \big)
    + \delta^{-1/3} \big( F^{5/6}N^{-1/2} + F^{-1/6}N^{1/2} \big),
  \]
  and the lemma follows on noting that, when $F \ll N^{3/2}$,
  \[
    F^{-1/6}N^{1/2} \ll F^{-1/3}N, \quad F^{5/6}N^{-1/2} \ll F^{1/6}N^{1/2}.
  \]
\end{proof}

Next, we turn to the bilinear sums needed in the proof of \eqref{2.5}. From now on, $X, X_1, N, H, J$ have the same
meaning as in \S\ref{s2} and $\varepsilon$ is subject to $0 < \varepsilon < \frac 12\big( \frac {16}{15} - c \big)$.

\begin{lemma}\label{l2}
  Suppose that $1 < c < \frac 65 - 6\varepsilon$, $M < M_1 \le 2M$, $2 \le K < K_1 \le 2K$, and
  \begin{equation}\label{3.7}
    M \ll X^{1 - 2c/3 - \varepsilon}.
  \end{equation}
  Further, suppose that $h, j$ are integers with $|h| \le H$, $|j| \le J$, $(h, j) \ne (0, 0)$, and that the
  coefficients $a_m$ satisfy $|a_m| \le 1$. Then
  \[
    \mathop{\sum_{M < m \le M_1} \sum_{K < k \le K_1}}_{X < mk \le X_1} a_m e\big( hm^ck^c + j(n - m^ck^c)^{\gamma} \big)
    \ll X^{2 - c - 4\varepsilon}.
  \]
\end{lemma}

\begin{proof}
  We shall focus on the case $j \ne 0$, the case $j = 0$ being similar and easier. We set
  \[
    y = jn^{\gamma}, \quad x = y^{-1}hn, \quad T = T_m = n^{\gamma}m^{-1} \asymp K.
  \]
  With this notation, we have
  \[
    f(k) = f_m(k) = hm^ck^c + j(n - m^ck^c)^{\gamma} = y \alpha( kT_m^{-1}),
  \]
  where
  \begin{equation}\label{3.1}
  \alpha(t) = \alpha(t; x) = xt^c + (1 - t^c)^{\gamma}.
  \end{equation}
  We have
  \begin{equation}\label{3.2}
    f''(k) = yT^{-2} \alpha''(kT^{-1}), \quad f'''(k) = yT^{-3} \alpha'''(kT^{-1}),
  \end{equation}
  and
  \begin{gather}
    \alpha''(t)  = (c - 1)t^{c - 2} \big( cx - (1 - t^c)^{\gamma - 2} \big), \label{3.2a}\\
    \alpha'''(t) = -(c - 1)(2c - 1)t^{2c - 3}(1 - t^c)^{\gamma - 3} + (c - 2)t^{-1}\alpha''(t). \label{3.2b}
  \end{gather}
  Moreover, by virtue of \eqref{2.1},
  \begin{equation}\label{3.8}
    \sfrac 12 < (kT^{-1})^c \le \sfrac 12(1.25)^c < \sfrac 45
  \end{equation}
  whenever $X < mk \le X_1$.

  Let $\delta_0 = X^{-\varepsilon/10}$. If $|x| \ge \delta_0^{-1}$, then by \eqref{3.2}, \eqref{3.2a}, and \eqref{3.8},
  \[
    |f''(k)| \asymp |xy|K^{-2} \asymp |h|nK^{-2} \qquad \implies \qquad JX^{1 - \varepsilon}K^{-2} \ll |f''(k)| \ll JXK^{-2}.
  \]
  Thus, by Lemma \ref{l1} with $\delta = X^{-\varepsilon}$, $F = JX$ and $N = K$,
  \begin{equation}\label{3.19}
    \mathop{\sum_{M < m \le M_1} \sum_{K < k \le K_1}}_{X < mk \le X_1} a_m e\big( f_m(k) \big)
    \ll MX^{\varepsilon/2} \big( X^{(c + \varepsilon)/6}K^{1/2} + KX^{-c/3} \big).
  \end{equation}
  Note that we need also to verify that $JX \le K^{3/2}$. This is a consequence of \eqref{3.7}.

  Suppose now that $|x| \le \delta_0^{-1}$. The set where $|\alpha''(kT^{-1})| \ge \delta_0$ consists of at most two intervals.
  Consequently, we can partition $[K, K_1]$ into at most three subintervals such that on each of them we have one of the
  following sets of conditions:
  \begin{itemize}
    \item  [i)] $\delta_0 |y|K^{-2} \ll |f''(k)| \ll \delta_0^{-1} |y|K^{-2}$;
    \item [ii)] $|y|K^{-3} \ll |f'''(k)| \ll |y|K^{-3}$, $|f''(k)| \ll \delta_0 |y|K^{-2}$.
  \end{itemize}
  Thus, by Lemma \ref{l1} with $\delta = \delta_0^2$, $F = \delta_0^{-1}|y| \asymp \delta_0^{-1}|j|X$, and $N = K$,
  \begin{equation}\label{3.20}
    \mathop{\sum_{M < m \le M_1} \sum_{K < k \le K_1}}_{X < mk \le X_1} a_m e\big( f_m(k) \big)
    \ll MX^{\varepsilon/10} \big( X^{(c + 2\varepsilon)/6}K^{1/2} + KX^{-1/3} \big).
  \end{equation}
  Again, we have $\delta_0^{-1}|j|X \le JX^{1 + \varepsilon/10} \le K^{3/2}$, by virtue of \eqref{3.7}.

  Combining \eqref{3.19} and \eqref{3.20}, we obtain the conclusion of the lemma, provided that $c < \frac 43 - 5\varepsilon$ and
  \[
    M \ll X^{3 - 7c/3 - 10\varepsilon}.
  \]
  Once again, the latter inequality is a consequence of \eqref{3.7}.
\end{proof}

\begin{lemma}\label{l3}
  Suppose that $1 < c < \frac {16}{15} - 2\varepsilon$, $M < M_1 \le 2M$, $K < K_1 \le 2K$, and
  \begin{equation}\label{3.9}
    X^{2c - 2 + 9\varepsilon} \ll M \ll X^{3 - 2c - 9\varepsilon}.
  \end{equation}
  Further, suppose that $h, j$ are integers with $|h| \le H$, $|j| \le J$, $(h, j) \ne (0, 0)$, and that the
  coefficients $a_m, b_k$ satisfy $|a_m| \le 1$, $|b_k| \le 1$. Then
  \[
    \mathop{\sum_{M < m \le M_1} \sum_{K < k \le K_1}}_{X < mk \le X_1} a_mb_k e\big( hm^ck^c + j(n - m^ck^c)^{\gamma} \big)
    \ll X^{2 - c - 4\varepsilon}.
  \]
\end{lemma}

\begin{proof}
  As in the proof of Lemma \ref{l2}, we shall focus on the case $j \ne 0$. By symmetry, we may assume that $M \ge X^{1/2}$. We set
  \[
    y = jn^{\gamma}, \quad x = y^{-1}hn, \quad T = n^{\gamma}.
  \]
  With this notation, we have
  \[
    f(k, m) = hm^ck^c + j(n - m^ck^c)^{\gamma} = y \alpha( mkT^{-1}),
  \]
  where $\alpha(t)$ is the function defined in \eqref{3.1}.

  By Cauchy's inequality and \cite[Lemma 8.17]{IwKo04},
  \begin{align}\label{3.10}
    \bigg| \mathop{\sum_{M < m \le M_1} \sum_{K < k \le K_1}}_{X < mk \le X_1} a_mb_k e\big( f(k, m) \big) \bigg|^2
    \ll \frac XQ \sum_{|q| \le Q} \sum_{K < k \le 2K} \bigg| \sum_{m \in \textbf{I}(k, q)} e\big( g(m; k, q) \big) \bigg| \notag \\
    \ll \frac {X^2}Q + \frac XQ \sum_{0 < |q| \le Q} \sum_{K < k \le 2K} \bigg| \sum_{m \in \textbf{I}(k, q)} e\big( g(m; k, q) \big) \bigg|,
  \end{align}
  where $g(m; k, q) = f(k + q, m) - f(k, m)$, $Q = J^2X^{6\varepsilon}$, and $\textbf{I}(k, q)$ is a subinterval of $[M, M_1]$ such that
  \[
    X < mk, m(k + q) \le X_1
  \]
  for all $m \in \textbf{I}(k, q)$. We remark that the right inequality in \eqref{3.9} ensures that $Q \ll KX^{-\varepsilon}$.
  When $q \ne 0$, we write
  \[
    g(m; k, q) = yT^{-1} \int_{mk}^{m(k + q)} \alpha'(tT^{-1}) \, dt =
    qy \int_0^1 \beta(m(k + \theta q)T^{-1}) \, \frac {d\theta}{k + \theta q},
  \]
  where $\beta(t) = t\alpha'(t)$. Introducing the notation
  \[
    z_{\theta} = z_{\theta}(k, q) = yq(k + \theta q)^{-1}, \quad U_{\theta} = U_{\theta}(k, q) = T(k + \theta q)^{-1} \asymp M,
  \]
  we find that
  \[
    g''(m) = \int_0^1 z_{\theta}U_{\theta}^{-2} \beta''(mU_{\theta}^{-1}) \, d\theta, \quad
    g'''(m) = \int_0^1 z_{\theta}U_{\theta}^{-3} \beta'''(mU_{\theta}^{-1}) \, d\theta,
  \]
  and
  \begin{gather}
    \beta''(t)  = (c - 1)t^{c - 2} \big( c^2x + (1 - t^c)^{\gamma - 3} (c + (c - 1)t^c) \big), \label{3.2c}\\
    \beta'''(t) = (c - 1)(2c - 1)t^{2c - 3}(1 - t^c)^{\gamma - 4}\big( (c - 1)t^c + 2c \big) + (c - 2)t^{-1}\beta''(t). \label{3.2d}
  \end{gather}

  Let $\delta_0 = X^{-\varepsilon/10}$. If $|x| \ge \delta_0^{-1}$, then by \eqref{3.2c} and a variant of \eqref{3.8},
  \[
    |g''(m)| \asymp |qxy|(XM)^{-1} \qquad \implies \qquad |q|JX^{-\varepsilon}M^{-1} \ll |g''(m)| \ll |q|JM^{-1}.
  \]
  Thus, by Lemma \ref{l1} with $\delta = X^{-\varepsilon}$, $F = |q|JM$ and $N = M$,
  \begin{equation}\label{3.23}
    \sum_{m \in \textbf{I}(k, q)} e\big( g(m; k, q) \big) \ll  (|q|J)^{1/6}M^{2/3}X^{\varepsilon/2}.
  \end{equation}
  Note that we need also to verify that $F \le M^{3/2}$, which holds if
  \begin{equation}\label{3.24}
    M \gg X^{6(c - 1) + 12\varepsilon}.
  \end{equation}

  Suppose now that $|x| \le \delta_0^{-1}$. We then deduce from \eqref{3.2c} and \eqref{3.2d} that
  \[
    |\beta''(mU_{\theta}^{-1})| \ll \delta_0^{-1}, \quad |\beta'''(mU_{\theta}^{-1})| \ll \delta_0^{-1},
  \]
  whence
  \[
    |\beta''(mU_{\theta}^{-1})| = |\beta''(mU_0^{-1})| + O \big( |q|K^{-1}\delta_0^{-1} \big) = |\beta''(mU_0^{-1})| + O \big( \delta_0^2 \big).
  \]
  We now note that the subset of $[M, M_1]$ where $|\beta''(mU_0^{-1})| \ge \delta_0$ consists of at most two intervals. Consequently,
  we can partition $[M, M_1]$ into at most three subintervals such that on each of them we have one of the following sets of conditions:
  \begin{itemize}
    \item  [i)] $\delta_0 |qy|(XM)^{-1} \ll |g''(m)| \ll \delta_0^{-1} |qy|(XM)^{-1}$;
    \item [ii)] $|qy|X^{-1}M^{-2} \ll |g'''(m)| \ll |qy|X^{-1}M^{-2}$, $|g''(m)| \ll \delta_0 |qy|(XM)^{-1}$.
  \end{itemize}
  Thus, Lemma \ref{l1} with $\delta = \delta_0^2$, $F = \delta_0^{-1}|qj|M$, and $N = M$ yields \eqref{3.23},
  provided that \eqref{3.24} holds.

  Combining \eqref{3.10} and \eqref{3.23}, we get
  \begin{equation}\label{3.25}
    \bigg| \mathop{\sum_{M < m \le M_1} \sum_{K < k \le K_1}}_{X < mk \le X_1} a_mb_k e\big( f(k, m) \big) \bigg|^2
    \ll X^2Q^{-1} + X^{2 + \varepsilon/2}(QJ)^{1/6}M^{-1/3}.
  \end{equation}
  In view of our choice of $Q$, the conclusion of the lemma follows from \eqref{3.25}, provided that
  \[
    M \gg X^{7.5(c - 1) + 10\varepsilon}.
  \]
  Both \eqref{3.24} and the last inequality follow from the assumption that $M \ge X^{1/2}$ and the hypothesis
  $c < \sfrac {16}{15} - 2\varepsilon$.
\end{proof}

We close this section with a lemma that will be needed in the proof of Theorem \ref{th2}.

\begin{lemma}\label{l4}
  Suppose that $1 < c < 2$, $2 \le X < X_1 \le 2X$, and $0 < \delta < \frac 14$. Let $\mathcal S_{\delta}$ denote the number
  of integers $n$ such that $X < n \le X_1$ and $\| n^c \| < \delta$. Then
  \[
    \mathcal S_{\delta} \ll \delta(X_1 - X) + \delta^{-1/2}X^{c/2}.
  \]
\end{lemma}

\begin{proof}
  Let $\Phi$ be the $1$-periodic extension of a smooth function that majorizes the characteristic function of the
  interval $[-\delta, \delta]$ and is majorized by the characteristic function of $[-2\delta, 2\delta]$. Then
  \begin{equation}\label{3.26}
    \mathcal S_{\delta} \le \sum_{X < n \le X_1} \Phi\big( n^c \big) = \sum_{X < n \le X_1} \hat\Phi(0)
    + \sum_{h \ne 0} \hat\Phi(h) \sum_{X < n \le X_1} e \big( hn^c \big).
  \end{equation}
  If $h \ne 0$, \cite[Corollary 8.13]{IwKo04} yields
  \[
    \sum_{X < n \le X_1} e \big( hn^c \big) \ll |h|^{1/2}X^{c/2},
  \]
  whence
  \begin{align}\label{3.27}
    \sum_{h \ne 0} \hat\Phi(h) \sum_{X < n \le X_1} e \big( hn^c \big) 
    &\ll X^{c/2} \sum_{h \ne 0} |\hat\Phi(h)||h|^{1/2} \notag\\
    &\ll X^{c/2} \sum_{h \ne 0} \frac {\delta|h|^{1/2}}{(1 + \delta|h|)^2} \ll \delta^{-1/2}X^{c/2}.
  \end{align}
  Since $\hat\Phi(0) \le 4\delta$, the lemma follows from \eqref{3.26} and \eqref{3.27}.
\end{proof}

\section{Proof of Theorem 1: conclusion}
\label{s4}

Suppose that $1 < c < \frac {16}{15}$ and $0 < \varepsilon < \frac 12\big( \frac {16}{15} - c \big)$. To prove \eqref{2.5},
we recall Vaughan's identity in the form of \cite[Proposition 13.4]{IwKo04}. We can use it to express the sum in
\eqref{2.5} as a linear combination of $O(\log^2X)$ sums of the form
\[
  \mathop{\sum_{M < m \le M_1} \sum_{K < k \le K_1}}_{X < mk \le X_1} a_mb_k e\big( hm^ck^c + j(n - m^ck^c)^{\gamma} \big),
\]
where either
\begin{itemize}
  \item  [i)] $|a_m| \ll m^{\varepsilon/2}$, $b_k = 1$, and $M \ll X^{2/3}$; or
  \item [ii)] $|a_m| \ll m^{\varepsilon/2}$, $|b_k| \ll k^{\varepsilon/2}$, and $X^{1/3} \ll M \ll X^{2/3}$.
\end{itemize}
A sum subject to conditions ii) is $\ll X^{2 - c - 3.5\varepsilon}$ by Lemma \ref{l3}. A sum subject to conditions i) can
be bounded using Lemma \ref{l2} if \eqref{3.7} holds and using Lemma \ref{l3} if \eqref{3.7} fails. In either case, the
resulting bound is $\ll X^{2 - c - 3.5\varepsilon}$. Therefore, each of the $O(\log^2X)$ terms in the decomposition of
\eqref{2.5} is $\ll X^{2 - c - 3.5\varepsilon}$. This establishes \eqref{2.5} and completes the proof of the theorem.

\section{Proof of Theorem \ref{th2}}
\label{s5}

We can cover the interval $( x^{1/2}, x ]$ by $O( (\log x)^3 )$ subintervals of the form $(N, N_1]$, with $N_1 =
N\big( 1 + (\log N)^{-2} \big)$. Thus, it suffices to show that
\begin{equation}\label{5.0}
  Z_c(N) \ll N^{3 - 3/c + 5\varepsilon/6},
\end{equation}
where $Z_c(N)$ is the number of integers $n$ in the range
\[
  N < n \le N \big( 1 + (\log N)^{-2} \big)
\]
that cannot be represented in the form \eqref{i.2}.

As in the proof of Theorem \ref{th1}, we derive solutions of \eqref{i.2} from solutions of \eqref{2.2}. We set
$\gamma = 1/c$, $\eta = (\log N)^{-2}$, and write
\[
  N_1 = (1 + \eta)N, \quad X = \big( \sfrac 12N \big)^{\gamma}, \quad  X_1 = (1 + \eta)X, \quad \delta = \gamma X^{1 - c}.
\]
Suppose that $N < n \le N_1$ and $X < p \le X_1$. Then
\[
  (1 - \eta)\delta < \gamma \big( n - p^c \big)^{\gamma - 1} < (1 + 2\eta)\delta.
\]
Assuming that $p$ satisfies the inequalities
\begin{equation}\label{4.1}
  4\eta < \{ p^c \} < 1 - 4\eta, \quad 1 - \delta - \eta\delta < \big\{ \big( n - p^c \big)^{\gamma} \big\} < 1 - \delta + \eta\delta,
\end{equation}
we deduce that
\begin{align*}
  \big( n - [p^c] \big)^{\gamma}
  &< \big( n - p^c \big)^{\gamma} + (1 - 4\eta)(1 + 2\eta)\delta + O \big( \delta n^{-1} \big) \\
  &< \big[ \big( n - p^c \big)^{\gamma} \big] + 1 - \eta\delta, \\
  \big( n + 1 - [p^c] \big)^{\gamma}
  &> \big( n - p^c \big)^{\gamma} + (1 + 4\eta)(1 - \eta)\delta + O \big( \delta n^{-1} \big) \\
  &> \big[ \big( n - p^c \big)^{\gamma} \big] + 1 + \eta\delta.
\end{align*}
In particular, a prime $p$, $X < p \le X_1$, that satisfies \eqref{4.1} yields a solution $m$ of \eqref{2.2} and
a representation of $n$ in the form \eqref{i.2}.

Let $\Phi$ be the $1$-periodic extension of a smooth function $\Phi_0$ that majorizes the characteristic function of $[6\eta, 1 - 6\eta]$
and is majorized by the characteristic function of $[4\eta, 1 - 4\eta]$. Further, let $\Psi$ be the $1$-periodic extension of
\[
  \Psi_0(t) = \psi_0\big( (2\eta\delta)^{-1}(t - 1 + \delta) + \sfrac 12),
\]
where $\psi_0$ is the function appearing in the proof of Theorem \ref{th1}. Then $\Psi_0$ is supported inside $[1 - \delta - \eta\delta,
1 - \delta + \eta\delta]$ and the Fourier coefficients of $\Psi$ satisfy
\begin{equation}\label{4.2}
  \hat\Psi(0) = 2\eta\delta, \quad |\hat\Psi(h)| \ll_r \eta\delta(1 + \eta\delta|h|)^{-r} \quad \text{for all } r \in \mathbb Z.
\end{equation}
Hence,
\begin{align}\label{4.3}
  \sum_{X < p \le X_1} \Phi\big( p^c \big) \Psi \big( (n - p^c)^{\gamma} \big) &=
  \sum_{h \in \mathbb Z} \sum_{X < p \le X_1} \Phi\big( p^c \big) \hat\Psi(h) \, e \big( h(n - p^c)^{\gamma} \big) \notag\\
  &= \hat\Psi(0) \sum_{X < p \le X_1} \Phi\big( p^c \big) + \mathcal R(n) \notag\\
  &= 2\eta\delta \big( \pi(X_1) - \pi(X) + O(\mathcal S) \big) + \mathcal R(n).
\end{align}
Here,
\[
  \mathcal R(n) = \sum_{h \ne 0} \hat\Psi(h) \sum_{X < p \le X_1} \Phi\big( p^c \big) \, e \big( h(n - p^c)^{\gamma} \big)
\]
and $\mathcal S$ is the number of integers $m$ such that $X < m \le X_1$ and $\| m^c \| < 6\eta$. By Lemma \ref{l4},
\begin{equation}\label{4.4}
  \mathcal S \ll \eta(X_1 - X) + \eta^{-1/2}X^{c/2} \ll \eta^2X.
\end{equation}
Combining \eqref{4.3}, \eqref{4.4} and the Prime Number Theorem, we find that
\begin{equation}\label{4.5}
  \sum_{X < p \le X_1} \Phi\big( p^c \big) \Psi \big( (n - p^c)^{\gamma} \big) \gg X^{2 - c}(\log X)^{-5}
\end{equation}
for any $n$, $N < n \le N_1$, for which we have
\begin{equation}\label{4.6}
  \mathcal R(n) \ll X^{2 - c - \varepsilon/12}.
\end{equation}
Since the sum on the right side of \eqref{4.5} is supported on the primes $p$ satisfying \eqref{4.1}, \eqref{5.0} will
follow if we show that \eqref{4.6} holds for all but $O\big( N^{3 - 3\gamma + 5\varepsilon/6} \big)$ integers $n \in (N, N_1]$.

Set $H = X^{c - 1 + \varepsilon/6}$. By \eqref{4.2} with $r = 2 + [2\varepsilon^{-1}]$, the contribution to $\mathcal R(n)$ from terms with $|h| > H$ is bounded. Consequently,
\[
  Z_c(N) \ll X^{-2 + \varepsilon/6} \sum_{N < n \le N_1} \mathcal R_1(n)^2,
\]
where
\[
  \mathcal R_1(n) = 
  \sum_{0 < |h| \le H} \bigg| \sum_{X < p \le X_1} \Phi \big( p^c \big) e \big( h(n - p^c)^{\gamma} \big) \bigg|.
\]
Appealing to Cauchy's inequality and the Weyl--van der Corput lemma \cite[Lemma 8.17]{IwKo04}, we obtain
\begin{align*}
  Z_c(N) &\ll X^{c - 3 + \varepsilon/3} \sum_{0 < |h| \le H} 
  \sum_{N < n \le N_1} \bigg| \sum_{X < p \le X_1} \Phi \big( p^c \big) e \big( h(n - p^c)^{\gamma} \big) \bigg|^2\\
  &\ll X^{c - 2 + \varepsilon/3}Q^{-1} \sum_{0 < |h| \le H} \sum_{|q| \le Q} \sum_{X < p \le X_1}
  \bigg| \sum_{N < n \le N_1} e(f(n)) \bigg|,
\end{align*} 
where $Q \le \eta X$ is a parameter at our disposal and
\[
  f(n) = qh\big( (n - p^c)^{\gamma} - (n - (p + q)^c)^{\gamma} \big). 
\]
We choose $Q = \eta X^{1 - \varepsilon/6}$. Then
\[
  |qh|N^{-1} \ll |f'(n)| \ll |qh|N^{-1} \ll \eta < \sfrac 12,
\]
so \cite[Corollary 8.11]{IwKo04} and the trivial bound yield
\[
  \sum_{N < n \le N_1} e(f(n)) \ll N(1 + |qh|)^{-1}.
\]
We conclude that
\begin{align*}
  Z_c(N) &\ll NX^{c - 2 + 2\varepsilon/3} \sum_{0 < |h| \le H} \sum_{|q| \le Q} (1 + |qh|)^{-1} 
  \ll NX^{2c - 3 + 5\varepsilon/6}.
\end{align*} 
This establishes \eqref{5.0} and completes the proof of the theorem.


\begin{thebibliography}{1}

\bibitem{BaGaZuX}
E.~P. Balanzario, M.~Z. Garaev, and R.~Zuazua, \emph{Exceptional set of a
  representation with fractional powers}, preprint.

\bibitem{Desh73}
J.-M. Deshouillers, \emph{Un probl\`eme binaire en th\'eorie additive}, Acta
  Arith. \textbf{25} (1973/74), 393--403.

\bibitem{Grit92}
S.~A. Gritsenko, \emph{Three additive problems}, Izv. Ross. Akad. Nauk
  \textbf{41} (1992), 447--464, in {R}ussian.

\bibitem{IwKo04}
H.~Iwaniec and E.~Kowalski, \emph{Analytic {Number} {Theory}}, American
  Mathematical Society, 2004.

\bibitem{Kony03}
S.~V. Konyagin, \emph{An additive problem with fractional powers}, Mat. Zametki
  \textbf{73} (2003), 633--636, in Russian.

\bibitem{Lapo99a}
M.~B.~S. Laporta, \emph{On a binary problem with prime numbers}, Math.
  Balkanica (N.S.) \textbf{13} (1999), 119--123.

\end{thebibliography}
\end{document}